\documentclass[11pt]{article}
\usepackage{amsmath,amssymb,amsthm,amscd,dsfont}
\usepackage{amsfonts,latexsym,rawfonts,amsmath,amssymb,amsthm}
\usepackage{lscape}
\usepackage{amscd, float,times,rotating}
\usepackage{pb-diagram}
\usepackage{hyperref}
\numberwithin{equation}{section}

\newtheorem{prop}{Proposition}[section]
\newtheorem{theorem}[prop]{Theorem}
\newtheorem{lemma}[prop]{Lemma}

\newtheorem{remark}[prop]{Remark}
\newtheorem{example}[prop]{Example}
\newtheorem{definition}[prop]{Definition}

\begin{document}
\title{Rigidity for inscribed radius estimate of asymptotically hyperbolic Einstein manifold}
\author{Xiaoshang Jin}
\date{}
\maketitle
\begin{abstract}
The inscribed radius of a compact manifold with boundary is bounded above if its Ricci curvature and mean curvature are bounded from below. The rigidity result implies that the upper bound can be achieved only in space form. In this paper, we generalize this result to asymptotically hyperbolic Einstein manifold. We get an upper bound of the relative volume of AH manifold and if we combine it with the recent work of Wang and Zhou, then the rigidity is obtained.
\end{abstract}
\section{Introduction}
Let $(N,g)$ be an $n+1$ dimensional complete manifold satisfying that $Ric_g\geq nk>0,$ the classic Myers theorem tells us that the diameter of $N$ satisfies $\mathrm{diam}(N,g)\leq\frac{\pi}{\sqrt{k}}.$ In 1975, Cheng \cite{cheng1975eigenvalue} showed that the equality holds if and only if $(N,g)$ is isometric to
$S^{n+1}_k.$ Here $S^{n+1}_k$ denote the $n+1$ dimensional standard sphere of curvature $k.$
\par If the lower bound of Ricci curvature is negative, then we could consider the manifold with boundary and have the following estimate of inscribed radius. Let $N$ be the $n+1$ dimensional manifold with smooth boundary $\partial N.$ If $Ric_g\geq -n$ and the mean curvature of $\partial N$ satisfies $H\geq nk>n$ for some constant $k>1,$  then
\begin{equation}\label{1.1}
  \mathrm{Inraj}(N,g)=\sup\limits_{x\in N} d_g(x,\partial N)\leq \coth^{-1} k.
\end{equation}
It is showed in \cite{kasue1983ricci} and \cite{li2015rigidity} that the equality holds if and only if $(N,g)$ is isometric to a geodesic ball of radius $\coth^{-1} k$ in hyperbolic space.
\par In this paper, we study the rigidity problem via inscribed radius estimate in AHE manifolds. Let's recall some definitions of AH and AHE manifolds.
\par \begin{definition}
Let $\overline{X}$ be an $n+1-$ dimensional compact manifold with interior $X$ and with boundary $\partial X=\Sigma$ which is a smooth $n-$ dimensional compact manifold. We say that a complete and noncompact manifold $(X,g^+)$ is a conformally compact manifold if there exists a smooth defining function $\rho$ on $\overline{X},$ such that
\begin{equation}\label{1.2}
  \rho>0\ \ in \ \ X,\ \ \ \ \rho=0\ \ on\ \ \Sigma,\ \ \ \ d\rho\neq 0\ \ on\ \ \Sigma.
\end{equation}
and  $g=\rho^2g^+$ can extend to a sooth Riemannian metric on $\overline{X}.$
\end{definition}
The induced metric $\hat{g}=g|_{\Sigma}$ is called the boundary metric. It is well-known that $g^+$ determines a conformal structure on the boundary $\Sigma$ and we call $(\Sigma,[\hat{g}])$ the conformal infinity of $g^+.$ If in addition $|d\rho|^2_g=1$ on the boundary $\Sigma,$ then we say that $(X,g^+)$ is an asymptotically hyperbolic (or AH for short) manifold.
\par Let $(X^{n+1},g^+)$ be a conformally compact manifold satisfying $Ric[g^+]=-ng^+,$ then $|d\rho|^2_g=1$ automatically holds on $\Sigma.$ We call $(X,g^+)$ an asymptotically hyperbolic Einstein (or AHE for short) manifold in this case.
\par For an AH manifold $(X,g^+)$ with conformal boundary $(\Sigma,\hat{g}),$ there exists a special defining function $x$ such that
the compactification $g=x^2g^+$ satisfies $|dx|^2_g=1$ in a neighbourhood of $\Sigma$ and we call $x$ the geodesic defining function.
\\
\par Since the AHE manifold is noncompact, we could consider the exhaust $\{E_x\}_{x\in(0,\delta)}$ of $X$ by compact sets. For each $E_x$ we define the difference of inscribed radius between $E_x$ and the geodesic ball in hyperbolic space. Then we study the limit as $x$ tends to $0$ and get the rigidity result. More concretely,
\begin{theorem}
  Suppose that $(X,g^+)$ is an $n+1-$ dimensional AHE manifold with minimized Yamabe boundary metric $(\Sigma,\hat{g})$ of positive scalar curvature and $x$ is the geodesic defining function. For sufficiently small $x>0,$ let
  \begin{equation}\label{1.3}
  F(x)=\coth^ {-1} {h_x}-\mathrm{Inraj}(E_x,g^+)
  \end{equation}
  where $E_x=X\setminus \Sigma\times(0,x)$ and $h_x=\min\limits_{p\in \partial E_x} \frac{H_x(p)}{n}$ and $H_x$ is the mean curvature of $\partial E_x.$ Then $F(x)\geq 0$ and the limit $\lim\limits_{x\rightarrow 0^+} F(x)=F(0)$ exists. Furthermore, $F(0)=0$ if and only if $(X,g^+)$ is isometric to the hyperbolic space.
\end{theorem}
Here is the outline of this paper. In section 2, we firstly review some materials about the relative volume of AH manifold, and then research some properties of the level sets in AH manifold, including the mean curvature, the inscribed radius and the cut locus. In section 3, we study the upper bound of the relative volume and prove a new relative volume inequality. As a consequence, the rigidity theorem can be obtained by utilizing the recent results of Wang and Zhou.
\section{The analysis on AH manifold}
In this section we assume that $(X,g^+)$ is an AH manifold satisfying
$$Ric[g^+]\geq -ng^+,\ \ R[g^+]+n(n+1)=o(x^2),$$
 and the boundary metric $(\Sigma,\hat{g})$ is the minimized Yamabe metric of positive scalar curvature. $g=x^2g^+$ is the geodesic comactification near conformal infinity. We could assume that $x\in[0,\delta)$ for some $\delta>0.$
\subsection{The relative volume function}
For any $p\in X,$ by the volume comparison theorem, the limit of $\frac{V(\partial B(p,t),g^+)}{\omega_n\sinh^n t}$ always exists as $t$ tends to infinity. In \cite{jin2022relative}, the author showed that the limit is dependent on $p$ for AH manifold. Then we could define the relative volume function $\mathcal{A}:X\rightarrow\mathbb{R}$ by
\begin{equation}\label{2.1}
   \mathcal{A}(p)=2^n\cdot\lim\limits_{t\rightarrow+\infty}e^{-nt}\cdot V(\partial B(p,t),g^+)
\end{equation}
Let $E$ be any compact set in $X,$ we can use the same method as lemma 4.1 in \cite{li2017gap} to show that the following limit also exists
 (one can see section 3 in \cite{jin2022relative} for more details).
 \begin{equation}\label{2.2}
   \mathcal{A}(E)=2^n\cdot\lim\limits_{t\rightarrow+\infty}e^{-nt}\cdot V(\partial B(E,t),g^+)
\end{equation}
Then we have the following estimate of relative volume (lemma 3.1 in \cite{jin2022relative}:
\begin{equation}\label{2.3}
  \forall p\in E,\mathcal{A}(p)\leq \mathcal{A}(E)\leq e^{n\cdot diam(E,g^+)}\mathcal{A}(p).
\end{equation}
We recall the important relative volume inequality in \cite{li2017gap}:
\begin{equation}\label{2.4}
\forall p\in X, \ (\frac{Y(\partial X,[\hat{g}])}{Y(\mathbb{S}^n,[g_{\mathbb{S}}])})^{\frac{n}{2}}\leq \frac{\mathcal{A}(p)}{\omega_n}\leq\frac{V(B(p,t),g^+)}{V(B(o,t),g^\mathbb{H})}\leq 1.
\end{equation}

\subsection{The mean curvature of the level set}
Set $X_x=\Sigma\times (0,x)$ for $x\in (0,\delta)$ and $E_x=X\setminus X_x.$ Let
$$\Sigma_x=\partial X_x=\partial E_x=\Sigma\times\{x\}$$
be the level set of the geodesic defining function.
Let $H_x$ denote the mean curvature of $\Sigma_x$ in $(E_x,g^+)$ with respect to the outer unit vector.
A direct calculation shows that
\begin{equation}\label{2.5}
  H_x=\Delta_{g^+}(-\ln x)=n-x\Delta_gx
\end{equation}
The standard formulas for conformal changes imply that the scalar curvature
\begin{equation}\label{2.6}
    S|_{\Sigma_x}=-2n\frac{\Delta_gx}{x}+tr_g(Ric[g^+]+ng^+)=-2n\frac{\Delta_gx}{x}+o(1).
\end{equation}

Hence
\begin{equation}
  H_x=n+\frac{x^2\cdot S|_{\Sigma_x}}{2n}+o(x^2)
\end{equation}
Set
\begin{equation}
s_x=\min\{S(p):p\in\Sigma_x\},\ \ \ \ \
  h_x=\min\limits_{p\in\Sigma_x}\frac{H_x(p)}{n}=1+\frac{x^2\cdot s_x}{2n^2}+o(x^2)
\end{equation}
When restricted on the boundary, by lemma 4.1 in \cite{jin2019boundary} (or lemma 1.2 in \cite{anderson2003boundary}),
\begin{equation}
s_0=S|_{\Sigma}=\frac{n}{n-1}S_{\hat{g}}.
\end{equation}
Recall we assume that the conformal infinity is of positive Yamabe type and we choose the minimized Yamabe metric $\hat{g},$ hence the scalar curvature $S_{\hat{g}}$ is a positive constant on $\Sigma.$ Then $H_x\geq nh_x>n$ at least for sufficiently small $x.$ Now we can consider $\coth^{-1}h_x$  and it denotes the radius of the geodesic ball in hyperbolic space whose mean curvature is $nh_x.$
Then
\begin{equation}
  \coth^{-1}h_x=\frac{1}{2}\ln (1+\frac{4n^2}{x^2\cdot s_x+o(x^2)})
\end{equation}
\begin{lemma}
  The function $\coth^{-1}h_x+\ln x$ is convergent to $\frac{1}{2}\ln\frac{4n(n-1)}{S_{\hat{g}}}$ as $x$ tends to $0.$
\end{lemma}
\begin{proof}
This is trivial as
 $$\lim\limits_{x\rightarrow 0}(\coth^{-1}h_x+\ln x)=\lim\limits_{x\rightarrow 0}\frac{1}{2}\ln (x^2+\frac{4n^2}{ s_x+o(1)})=\frac{1}{2}\ln\frac{4n(n-1)}{S_{\hat{g}}}.$$
\end{proof}

\subsection{The inscribed radius and the cut locus of level sets}
Since $E_x$ is closed of the compact manifold $\overline{X},$ $E_x$ and its boundary $\Sigma_x$ are also compact sets. Consider the inscribed radius
of $E_x:$
\begin{equation}
\mathrm{Inraj}(E_x,g^+)=\sup\limits_{p\in E_x}d_{g^+}(p,\Sigma_x).
\end{equation}
\begin{definition}
Foot point: For any $p\in E_x,$ by the compactness of $(E_x,g^+),$ there exists at least one point $q\in \Sigma_x$ such that
$d_{g^+}(p,\Sigma_x)=d_{g^+}(p,q).$ We call $q$ a foot point on $\Sigma_x$ of $p.$
\end{definition}
For any $q\in \Sigma_x,$ we set $\sigma_q:[0,T)\rightarrow E_x$ to be the normal geodesic satisfying $\sigma_q(0)=q$ and $\sigma_q'(0)\bot T_q\Sigma_x.$
If $q$ is a foot point on $\Sigma_x$ of some point $p\in E_x,$ then $\sigma_q:[0,l]\rightarrow E_x$ is the unique normal minimal geodesic connecting $q$ and $p$ where $l=d_{g^+}(p,q).$
\par We define the function $\tau_x:\Sigma_x\rightarrow \mathbb{R}$ by
\begin{equation}
\tau_x(q)=\sup\{t>0:d_{g^+}(\sigma_q(t),\Sigma_x)=t\}.
\end{equation}
Then it is easy to prove that $\mathrm{Inraj}(E_x,g^+)=\sup\{\tau_x(q):q\in \Sigma_x\}.$
\begin{definition}
Cut locus: We call $Cut(\Sigma_x)=\{\sigma_q(\tau_x(q)):q\in\Sigma_x\}$ the cut
locus for the boundary $\Sigma_x$ in $E_x.$
\end{definition}
The classic Riemann geometry theory tells us that $Cut(\Sigma_x)$ is a closed subset and the distance function $d_{g^+}(\Sigma_x,\cdot)$ is smooth in $E_x\setminus Cut(\Sigma_x).$ Furthermore, if $p=\sigma_q(\tau_x(q))\in Cut(\Sigma_x),$ then either $p$ is the first focal point of $\Sigma_x$ along $\sigma_q,$ or there exists another foot point $q'$ on $\Sigma_x$ of $p.$

\begin{lemma}
  For any different number $x_1, x_2\in (0,\delta),$ the cut locus are the same: $Cut(\Sigma_{x_1})=Cut(\Sigma_{x_2}).$
\end{lemma}
\begin{proof}
  We may assume that $x_1<x_2.$ For any $q_1\in\Sigma_{x_1},$ consider the normal geodesic $\sigma_{q_1}(t)$ with $\sigma_{q_1}(0)=q_1$ and $\sigma_{q_1}'(0)\bot T_{q_1}\Sigma_{x_1}.$ Notice that $g(\nabla x,\nabla x)=1$ on $\overline{X}_\delta$ implies that
\begin{equation}
g^+(\nabla[g^+](-\ln x),\nabla[g^+](-\ln x))=1.
\end{equation}
 Hence $-\ln x$ is the distance function on $(X_\delta,g^+).$ $\Sigma_{x_1}$ is the level set of the distance function $x$ in $\overline{X}_\delta,$ and it is also the level set of $-\ln x$ in $E_{x_1}.$ We get that
  $\sigma_{q_1}'(0)= \nabla[g^+](\ln x)|_{q_1}$ and $d_{g^+}(q_1,\Sigma_{x_2})=\ln\frac{x_2}{x_1}.$ By gauss lemma,
\begin{equation}
\sigma_{q_1}(\ln\frac{x_2}{x_1})=q_2\in \Sigma_{x_2},\ \ \ \sigma_{q_1}'(\ln\frac{x_2}{x_1})=\nabla[g^+](\ln x)|_{q_2}\bot T_{q_2}\Sigma_{x_2}.
\end{equation}
  Then $\tau_{x_1}(q_1)=\tau_{x_2}(q_2)+\ln\frac{x_2}{x_1},$ which would imply that $\sigma_{q_1}(\tau_{x_1}(q_1))=\sigma_{q_2}(\tau_{x_2}(q_2))$ Hence $Cut(\Sigma_{x_1})\subseteq Cut(\Sigma_{x_2}).$ The opposite direction could be obtained in the same way.
\end{proof}
From the above lemma, we may define the cut locus of $(\Sigma,\hat{g})$ in $X$ by
\begin{equation}
Cut(\Sigma)=Cut(\Sigma_x)
\end{equation}
 for any $x\in (0,\delta).$ We also get that for any $p\in E_{x_2}\subseteq E_{x_1},$
\begin{equation}
d_{g^+}(p,\Sigma_{x_1})=d_{g^+}(p,\Sigma_{x_2})+\ln\frac{x_2}{x_1}.
\end{equation}
Then $d_{g^+}(p,\Sigma_x)+\ln x=c(p)$ is a constant only depending on $p.$ As a consequence, we deduce the following lemma:
\begin{lemma}
$\mathrm{Inraj}(E_x,g^+)+\ln x$ is a constant function for $x\in (0,\delta).$
\end{lemma}
Furthermore, the constant $c_0=\mathrm{Inraj}(E_x,g^+)+\ln x=\max\limits_{p\in Cut(\Sigma)} c(p).$
Combining lemma 2.1 and lemma 2.5, we get that
\begin{equation}
\begin{aligned}
  \lim\limits_{x\rightarrow 0^+}F(x)&=\lim\limits_{x\rightarrow 0^+}[(\coth^{-1}h_x+\ln x)-(\mathrm{Inraj}(E_x,g^+)+\ln x)]
                                  \\&=\frac{1}{2}\ln\frac{4n(n-1)}{S_{\hat{g}}}-c_0
  \end{aligned}
\end{equation}
Now we claim that $F(0)$ is a scaling invariant. In fact, let $\hat{h}=K^2\hat{g}\in [\hat{g}]$ for some positive constant $K.$
Then $Kx$ is the unique geodesic defining function of $X$ with boundary metric $(\Sigma,\hat{h})$ and the scalar curvature $S_{\hat{h}}=K^{-2}S_{\hat{g}}.$ Let $h=(Kx)^2g^+=K^2g,$ then $\Sigma\times (0,\varepsilon)_g=\Sigma\times (0,K\varepsilon)_h$ for small $\varepsilon>0.$ If we use $E_y(h)$ to denote the complementary set of $\Sigma\times (0,y)_h$ in $X,$ then the constant
\begin{equation}
c_0(\hat{h})\equiv \mathrm{Inraj}(E_y(h),g^+)+\ln y=\mathrm{Inraj}(E_{\frac{y}{K}}(g),g^+)+\ln y=c_0(\hat{g})+\ln K.
\end{equation}
Hence
\begin{equation}
F_{\hat{h}}(0)=\frac{1}{2}\ln\frac{4n(n-1)}{S_{\hat{h}}}-c_0(\hat{h})=F_{\hat{g}}(0)
\end{equation}
Then we obtained that $F(0)$ is a scaling invariant.
\subsection{Some examples}
\begin{example}
  Hyperbolic space: (Poincar\'e ball model)
  $$(\mathbb{B}^{n+1},\ g_\mathbb{H}=x^{-2}(dx^2+\frac{(4-x^2)^2}{16}g_{\mathbb{S}^n}),$$
  Here $x$ is the geodesic defining function and the conformal infinity is the standard conformal round sphere $(\mathbb{S}^n,[g_{\mathbb{S}^n}]).$
\end{example}
We choose the boundary Yamabe metric $g_{\mathbb{S}^n}.$ Then $x$ is defined on $[0,2)$ and $Cut(\mathbb{S}^n)=\{o\}$ is a single point.
If we let $x\rightarrow 2,$ then
$$c_0=\mathrm{Inraj}(E_x,g_\mathbb{H})+\ln x=\ln 2.$$ Hence
$$F(0)=\frac{1}{2}\ln\frac{4n(n-1)}{S(g_{\mathbb{S}^n})}-c_0=0. $$

\begin{example}
  Hyperbolic manifold:

  $$(\mathbb{R}^n\times \mathbb{S}^1(\lambda),g^+=dt^2+\sinh^2tg_{\mathbb{S}^{n-1}}+\cosh^2tg_{\mathbb{S}^1(\lambda)}).
$$
Here $n\geq 3,$ and we use the polar coordinate system $\mathbb{R}^n=[0,+\infty)\times\mathbb{S}^{n-1}.$
The conformal infinity of $g^+$ is
$$(\mathbb{S}^{n-1}\times\mathbb{S}^1,[g_{\mathbb{S}^{n-1}}+g_{\mathbb{S}^1(\lambda)}]).$$
Let $2x=e^{-t},$ then
$g^+=x^{-2}[dx^2+(x^2-1)^2g_{\mathbb{S}^{n-1}}+(x^2+1)^2g_{\mathbb{S}^1(\lambda)}].$
$x$ is the geodesic defining function on $[0,1)$ and the boundary metric is $\hat{g}=g_{\mathbb{S}^{n-1}}+g_{\mathbb{S}^1(\lambda)}.$
\end{example}
If $\lambda\in(0,(n-2)^{-\frac{1}{2}}2\pi],$ then $\hat{g}$ is the unique Yamabe metric (\cite{schoen1989variational}). Then
 $$Cut(\mathbb{S}^{n-1}\times\mathbb{S}^1)=\{0_n\}\times\mathbb{S}^1(2\lambda)$$ and $0_n$ is the origin of $\mathbb{R}^n.$
 Let $x\rightarrow 1,$ we get
$$c_0=\mathrm{Inraj}(E_x,g^+)+\ln x=0$$ and hence
$$F(0)=\frac{1}{2}\ln\frac{4n(n-1)}{S_{\hat{g}}}-c_0=\frac{1}{2}\ln\frac{4n}{n-2}$$
doesn't depend on $\lambda.$
\section{Proof of theorem 1}

\subsection{A new relative volume inequality}
\begin{theorem}
  Let $(X,g^+)$ be an $n+1-$ dimensional AH manifold with positive boundary minimized Yamabe metric $(\Sigma,\hat{g}).$ Assume that $x$
  is the geodesic defining function and
  $$Ric[g^+]\geq ng^+,\ \ R[g^+]+n(n+1)=o(x^2).$$
  Let $C=Cut(\Sigma)$ be the cut locus of $(\Sigma,\hat{g}).$ Then there exists  $p\in C$ such that
  \begin{equation}
   \frac{\mathcal{A}(p)}{\omega_n}\leq(\frac{Y(\Sigma,[\hat{g}])}{Y(\mathbb{S}^n,[g_{\mathbb{S}}])})^{\frac{n}{2}}\cdot e^{nF(0)}\leq \frac{\mathcal{A}(C)}{\omega_n}.
  \end{equation}
\end{theorem}
We use $B(p,t)$ (resp. $B(C,t)$) to denote the geodesic ball of radius $t$ centred at $p$ (resp.$C$) in $(X,g^+)$ and $V$ to denote the volume with respect to $g^+.$ Let $g=x^2g^+$ be the geodesic compactification and $s_C(\cdot)=d_{g^+}(C,\cdot)$ be the distance function of $C.$  Set $g^C=e^{-2s_C}g^+.$ Then
\begin{equation}
  \begin{aligned}
     \mathcal{A}(C)&=2^n\lim\limits_{t\rightarrow\infty}e^{-nt}V(\partial B(C,t))=2^n\int_\Sigma dV_{g^C}
     \\&=2^n\lim\limits_{x\rightarrow 0}\int_{\Sigma_x}dV_{g^C}=2^n\lim\limits_{x\rightarrow 0}\int_{\Sigma_x}\frac{1}{(xe^{s_C})^n}dV_g
     \\&\geq 2^n\lim\limits_{x\rightarrow 0}\int_{\Sigma_x}\frac{1}{(xe^{\mathrm{Inraj}(E_x)})^n}dV_g
     \\&= 2^n\lim\limits_{x\rightarrow 0}\int_{\Sigma_x}\frac{1}{(xe^{\coth^{-1}h_x})^n}\cdot e^{n[\coth^{-1}h_x-{\mathrm{Inraj}(E_x)}]}dV_g
     \\&=2^n\int_\Sigma(\frac{S_{\hat{g}}}{4n(n-1)})^{\frac{n}{2}}\cdot e^{nF(0)} dV_{\hat{g}}
     =\omega_n\cdot(\frac{Y(\Sigma,[\hat{g}])}{Y(\mathbb{S}^n,[g_{\mathbb{S}}])})^{\frac{n}{2}}\cdot e^{nF(0)}
  \end{aligned}
\end{equation}
 On the other hand, firstly we recall the Heintze Karcher inequality (theorem 2.1 in \cite{heintze1978general}),
 \begin{equation}\label{3.3}
   \frac{V(E_x)}{V(\Sigma_x)}\leq\int_0^{\coth^{-1}h_x}(\cosh r-h_x\sinh r)^ndr
 \end{equation}
 Where $h_x=\min\frac{H|_{\Sigma_x}}{n}.$ Let $x$ tends to 0, then
 \begin{equation}\label{3.4}
 \begin{aligned}
   \lim\limits_{x\rightarrow 0} \int_0^{\coth^{-1}h_x} &(\cosh r-h_x\sinh r)^ndr
   = \lim\limits_{y\rightarrow+\infty}\int_0^y(\cosh r-\coth y\sinh r)^ndr
     \\&= \lim\limits_{y\rightarrow+\infty}\frac{\int_0^y \sinh^n(y-r)dr}{\sinh^ny}
     =\lim\limits_{y\rightarrow+\infty}\frac{\int_0^y \sinh^n tdt}{\sinh^ny}=\frac{1}{n}.
   \end{aligned}
 \end{equation}
 For the left side, there exists a point $p\in C$ such that for all $x\in(0,\delta),\ d_{g^+}(p,\Sigma_x)={\mathrm{Inraj}(E_x)}.$ Then $E_x\supseteq B(p,{\mathrm{Inraj}(E_x)}).$ Notice that $V(\Sigma_x)=x^{-n}V(\Sigma_x,g),$ we deduce

  \begin{equation}\label{3.5}
 \begin{aligned}
  \varlimsup \limits_{x\rightarrow 0}\frac{V(E_x)}{V(\Sigma_x)} & \geq  \varlimsup_{x\rightarrow 0}\frac{V[B(p,{\mathrm{Inraj}(E_x)})]}{x^{-n}V(\Sigma_x,g)}
 \\ &=   \varlimsup_{x\rightarrow 0}\frac{V[B(p,{\mathrm{Inraj}(E_x)})]}{e^{n{\mathrm{Inraj}(E_x)}}}\cdot\frac{1}{V(\Sigma_x,g)}\cdot(xe^{{\mathrm{Inraj}(E_x)}})^n
 \\&=\frac{\mathcal{A}(p)}{n2^n}\cdot \frac{1}{V(\Sigma,\hat{g})}\cdot(\frac{4n(n-1)}{S_{\hat{g}}})^{\frac{n}{2}}e^{-nF(0)}
 \\&=\frac{1}{n}\frac{\mathcal{A}(p)}{\omega_n}\cdot(\frac{Y(\Sigma,[\hat{g}])}{Y(\mathbb{S}^n,[g_{\mathbb{S}}])})^{-\frac{n}{2}}\cdot e^{-nF(0)}
   \end{aligned}
 \end{equation}
  Then we finish the proof of theorem 3.1 from (\ref{3.3}),(\ref{3.4}) and (\ref{3.5}).

 \begin{remark}
   If we use $diam(C,g^+)=\sup \{d_{g^+}(p,q):p,q\in C\}$ to denote the diameter of $C,$ then from theorem 3.1 and (\ref{2.3})
   we get that
 \begin{equation}\label{3.6}
   e^{n[F(0)-diam(C,g^+)]}(\frac{Y(\Sigma,[\hat{g}])}{Y(\mathbb{S}^n,[g_{\mathbb{S}}])})^{\frac{n}{2}}\leq \frac{\mathcal{A}(p)}{\omega_n}
   \leq e^{nF(0)}(\frac{Y(\Sigma,[\hat{g}])}{Y(\mathbb{S}^n,[g_{\mathbb{S}}])})^{\frac{n}{2}}.
  \end{equation}
  On one hand, (\ref{3.6}) provides an upper bound of the relative volume. On the other hand, the lower bound of $\frac{\mathcal{A}(p)}{\omega_n}$ in (\ref{3.6}) is better than (\ref{2.4}) sometimes. If we consider example 2.7 of the hyperbolic manifold, then $F(0)=\frac{1}{2}\ln\frac{4n}{n-2}$ and $diam(C,g^+)=2\pi\lambda.$ Hence $F(0)>diam(C,g^+)$ as long as $\lambda$ is small.
 \end{remark}

\subsection{The fractional Yamabe problem on AH manifold}
In order to prove the rigidity result of theprem, here we would introduce some background materials about the fractional Yamabe constant.  Let $(X^{n+1},g^+)$ be an AH manifold with boundary $(\Sigma,\hat{g})$ and $x$ is the unique geodesic defining function near infinity. It is showed in \cite{mazzeo1987meromorphic} that given $f\in C^\infty(\Sigma),$ and $s\in \mathbb{C}$ satisfying $Re(s)>\frac{n}{2}$ and $s(n-s)$ is not an $L^2$-eigenvalue for $-\Delta_{g^+},$ the Poisson equation
$$
  -\Delta_{g^+}u-s(n-s)u=0\ \ in\ X
$$
has a unique solution of the form
$$
  u= Fx^{n-s}+Gx^s, \ \ F,G\in C^\infty(\overline{X}), \ \ F|_\Sigma=f
$$
The scattering operator on $\Sigma$ is then defined as
$$
S(s)f=G|_\Sigma.
$$
Let $\gamma\in (0,\frac{n}{2}),$ Graham and Zworski defined the fractional GJMS operator
$P_{\gamma}^{\hat{g}}$ as the normalized scattering operator in \cite{graham2003scattering}
\begin{equation}
P_{\gamma}^{\hat{g}}f=d_\gamma P(\frac{n}{2}+\gamma)f,\ \   d_\gamma=2^{2\gamma}\frac{\Gamma(\gamma)}{\Gamma(-\gamma)}
\end{equation}
It is a conformally covariant operator and we call
\begin{equation}
Q_{\gamma}^{\hat{g}}=\frac{2}{n-2\gamma}P_{\gamma}^{\hat{g}}1
\end{equation}
the fractional scalar curvature. Then the fractional Yamabe constant is defined as:
\begin{equation}
\Lambda_{\gamma}(\Sigma,[{\hat{g}}])=\inf\limits_{f\in C^\infty(\Sigma),f>0}\frac{\int_\Sigma fP_{\gamma}^{\hat{g}}fdV_{\hat{g}}}{(\int_\Sigma f^{\frac{2n}{n-2\gamma}}dV_{\hat{g}})^{\frac{n-2\gamma}{n}}}=\inf\limits_{\hat{h}\in [\hat{g}]}\frac{\frac{n-2\gamma}{2}\int_\Sigma Q_{\gamma}^{\hat{h}}dV_{\hat{h}}}{(\int_\Sigma dV_{\hat{h}})^{\frac{n-2\gamma}{n}}}
\end{equation}
The
 $\gamma-$Yamabe problem is to find a metric in the conformal class $[\hat{g}]$ that $\Lambda_{\gamma}(\Sigma,[{\hat{g}}])$ is achieved.
It is clear that the solution $\hat{g}$ has a constant fractional
scalar curvature.
If $\gamma=1,$ then $P_{\gamma}^{\hat{g}}$ is the conformal Laplacian and
\begin{equation}\label{3.10}
\Lambda_{1}(\Sigma,[{\hat{g}}])=\frac{n-2}{4(n-1)}Y(\Sigma,[{\hat{g}}]).
\end{equation}
If $\gamma=\frac{1}{2},$  then
\begin{equation}
\Lambda_{\frac{1}{2}}(\Sigma,[{\hat{g}}])=\frac{n-1}{4n}Q(\overline{X},\Sigma;[g]).
\end{equation}
Here
\begin{equation}
Q(\overline{X},\Sigma;[g]):=\inf\limits_{h\in[g]} \frac{\int_{\overline{X}} R_h dV_h+2\int_\Sigma H_hdV_{h|_\Sigma}}{(Vol(\Sigma,h|_\Sigma) )^{\frac{n-1}{n}}}
\end{equation}
is the second type of Yamabe constant of compact manifold with boundary and was introduced by Escobar in \cite{escobar1992conformal} .
The readers could find the reference \cite{case2016fractional}\cite{gonzalez2013fractional}\cite{joshi2000inverse} for more information about the fractional GJMS operator and $\gamma-$Yamabe problem.
\par For $\gamma\in (0,1),$ we refer the recent work about the fractional Yamabe constant by Wang and Zhou.
\begin{lemma}[theorem 1.2 in \cite{wang2021lower}]
Assume that $(X,g^+)$ is an $n+1-$ dimensiona AHE metric with boundary metric $(\Sigma,\hat{g})$ of positive Yamabe type, let $p\in X$ and $\gamma\in (0,1),$ then
\begin{equation}\label{3.13}
(\frac{\Lambda_\gamma(\Sigma,[\hat{g}])}{\Lambda_\gamma(\mathbb{S}^n,[g_{\mathbb{S}}])})^{\frac{n}{2\gamma}}
\leq\frac{\mathcal{A}(p)}{\omega_n}\leq 1.
\end{equation}
\end{lemma}

\begin{lemma}[theorem5 in \cite{chen2019escobar},theorem 1.1 in \cite{wang2021comparison}]
Suppose $(X,g^+)$ is an $n+1-$ dimensional $C^{3,\alpha}$ AHE metric with boundary metric $(\Sigma,\hat{g}).$
Assume $\frac{n^2-1}{4}\notin Spec(-\Delta_{g^+})$ and the $\frac{1}{2}-$Yamabe problem is solvable for $(\Sigma,[\hat{g}]),$ then
\begin{equation}\label{3.14}
\frac{\Lambda_1(\Sigma,[\hat{g}])}{\Lambda_1(\mathbb{S}^n,[g_{\mathbb{S}}])}\leq
(\frac{\Lambda_{\frac{1}{2}}(\Sigma,[\hat{g}])}{\Lambda_{\frac{1}{2}}(\mathbb{S}^n,[g_{\mathbb{S}}])})^2
\end{equation}
and the equality holds if and only if $(X, g^+)$ is isometric to the hyperbolic space.
\end{lemma}

\subsection{The rigidity of AHE manifold}
We provide the proof for the rigidity result of theorem 1.2.
\par Case 1: $\Lambda_{\frac{1}{2}}(\Sigma,[\hat{g}])=\Lambda_{\frac{1}{2}}(\mathbb{S}^n,[g_{\mathbb{S}}]).$ Then by lemma 3.3, $\mathcal{A}(p)=\omega_n.$
The Bishop-Gromov volume comparison theorem implies that $X$ is isometric to a hyperbolic space (theorem 4.5 in \cite{jin2022relative}).
\par Case 2: $\Lambda_{\frac{1}{2}}(\Sigma,[\hat{g}])<\Lambda_{\frac{1}{2}}(\mathbb{S}^n,[g_{\mathbb{S}}]).$ Then the $\frac{1}{2}-$Yamabe problem is solvable for $(\Sigma,[\hat{g}]),$ (see proposition 2.1 in \cite{escobar1992conformal} or theorem 1.4 in \cite{gonzalez2013fractional}). We also have that
$Spec(-\Delta_{g^+})$ has no $L^2$ eigenvalue provided $Y(\Sigma,[\hat{g}])\geq 0$ by \cite{lee1994spectrum}. Hence (\ref{3.14}) holds.
Then (\ref{3.6}),(\ref{3.10}),(\ref{3.13}) and (\ref{3.14}) tell us that $\exists p\in Cut(\Sigma),$
\begin{equation}
  (\frac{Y(\Sigma,[\hat{g}])}{Y(\mathbb{S}^n,[g_{\mathbb{S}}])})^{\frac{n}{2}}\leq
(\frac{\Lambda_{\frac{1}{2}}(\Sigma,[\hat{g}])}{\Lambda_{\frac{1}{2}}(\mathbb{S}^n,[g_{\mathbb{S}}])})^n
\leq\frac{\mathcal{A}(p)}{\omega_n}\leq e^{nF(0)}(\frac{Y(\Sigma,[\hat{g}])}{Y(\mathbb{S}^n,[g_{\mathbb{S}}])})^{\frac{n}{2}}.
\end{equation}
If $F(0)=0,$ then the equality is achieved and $(X, g^+)$ is isometric to the hyperbolic space by lemma 3.4.
\bibliographystyle{plain}%

\bibliography{bibfile}

\noindent{Xiaoshang Jin}\\
  School of mathematics and statistics, Huazhong University of science and technology, Wuhan, P.R. China. 430074
 \\Email address: jinxs@hust.edu.cn

\end{document}